# Planar Analytic Functions

L. Gerritzen

(2.1.2007)

### Abstract

If $a$ is a point in the domain of convergence of a planar power series $f$ in a single variable $x$ one con expand $f$ into a planar power series in the variable $(x-a)$. One arrives at the notion of planar analytic functions on any domain D in the complex plane. It can be described by section $S$ of the sheaf of planar germs. The $k$ – ary exponential series $\exp_k(x)$ has infinite radius of convergence. It is possible to define a planar analogue $\wp_\mathbb{P}$ of the classical zeta-function. As yet a functional equation for $\wp_\mathbb{P}$ has not been obtained.

# 1   Radius of convergence

Let $\mathbb{C}\{\{x\}\}_\mathbb{P}$ be the $\mathbb{C}$ - algebra of formal planar power series in a variable x over the field $\mathbb{C}$ of complex numbers, see [G1], An element $f \in \mathbb{C}\{\{x\}\}_\mathbb{P}$ has a unique expansion $f = \sum_{T \in \mathbb{P}} \gamma_T x^T$ where $\mathbb{P}$ is the set of finite, planar, reduced rooted trees and $\gamma_T \in \mathbb{C}$ for all T. The coefficient $\gamma_T$ of is also denoted by $\langle f, x^T \rangle$ and thus $f = \sum_{T \in \mathbb{P}} (f, x^T) x^T$.

In this way $\mathbb{C}\{x\}\}_\mathbb{P}$ is identified with $\mathbb{C}$ - vector space of $\mathbb{C}$ - valued functions on $\mathbb{P}$ or on $\{x^T : T \in \mathbb{P}\}$.

Definition:   $\operatorname{rad}(f) := \sup \{r \in \mathbb{R}_{\geq 0} : \sum_{T \in \mathbb{P}} |\gamma_T| r^{\deg(T)} < \infty \}$ where $\deg(T)$ is the

number of leaves of T and $|\gamma_T|$ is the absolute value of $\gamma_T$. Then $\operatorname{rad}(f) \in \mathbb{R}_{\geq 0} \cup \{\infty\}$ is called the radius of convergence of $f$.

Let now $(f_n)_{n \geq 0}$ be a sequence of planar power series $f_n \in \mathbb{C}\{\{x\}\}_\mathbb{P}$. It is called point wise convergent, if $\lim \langle T, f_n \rangle$ for $n \to \infty$
is convergent in $\mathbb{C}$ for all $T \in \mathbb{P}$.

In this case the power series
$$\lim_{n\to\infty} f_n := \sum_{T\in\mathbb{P}} \lim_{n\to\infty} \langle T, f_n\rangle \cdot x^T$$
is the limit of the sequence $(f_n)_{n\geq 0}$.

## 2 Expansion into powers of (x – a)

Let $\mathbb{P}_n$ be the set of trees T in $\mathbb{P}$ with $\deg(T) \leq n$.
Let $\mathbb{C}\{x\}_{\leq n}$ be the $\mathbb{C}$ - vector pace of planar polynomials in x of degree $\leq n$. Then by definition the system $\{x^T : T\in\mathbb{P}_n\}$ is a $\mathbb{C}$-basis of $\mathbb{C}\{x\}_{\leq n}$.

**Proposition 2.1.** Let $a \in \mathbb{C}$. Then $\{(x-a)^T : T\in\mathbb{P}_n\}$ is a $\mathbb{C}$ - basis of $\mathbb{C}\{x\}_{\leq n}$. The coefficient $\langle f, (x-a)^T\rangle$ of $f \in \mathbb{C}\{x\}_{\leq n}$ relative to this basis and to the basis element $(x-a)^T$ is given by the formula
$$\langle f, (x-a)^T\rangle = \sum_{\in \mathbb{P}_n} \langle f, x^n\rangle \cdot (U/T)\, a^{(\deg(U) - \deg(T))}$$
where (U/T) is the planar binominal coefficient of U over T, see [G7],

Denote by $\mathbb{C}\{\{x-a\}\}_{\mathbb{P}}$ the $\mathbb{C}$-algebra of formal planar power series in the variable x-a.

An element $f \in \mathbb{C}\{\{x-a\}\}_{\mathbb{P}}$ has a unique expansion $f = \sum_{T\in\mathbb{P}} \gamma_T (x-a)^T$
with $\gamma_T \in \mathbb{C}$. The coefficient $\gamma_T$ is also denoted by $\langle f, (x-a)^T\rangle$ and thus
$$f = \sum_{T\in\mathbb{P}} \langle f, (x-a)^T\rangle \cdot (x-a)^T$$

In this way the set of $\mathbb{C}\{\{x-a\}\}_{\mathbb{P}}$ is identified with the set of $\mathbb{C}$ - valued functions on $\{(x-a)^T : T\in\mathbb{P}\}$.

A sequence $(f_n)_{n\geq 0}$ in $\mathbb{C}\{\{x-a\}\}_{\mathbb{P}}$ is pointwise convergent,
if $\lim_{n\to\infty} \langle f_n, (x-a)^T\rangle$ exists in $\mathbb{C}$ for all $T\in\mathbb{P}$.

Let $f \in \mathbb{C}\{\{x\}\}_{\mathbb{P}}$ and $a \in \mathbb{C}$.

Let $f_n := \sum_{T\in\mathbb{P}_n} \langle f, x^T\rangle \cdot x^T$ and

$$g_n := \sum_{T\in\mathbb{P}_n} \langle f_n, (x-a)^T\rangle (x-a)^T$$

Then $g_n = f_n$ and $\lim_{n\to\infty} f_n = f$ in $\mathbb{C}\{\{x\}\}_{\mathbb{P}}$.

Now we consider the sequence $(g_n)_{n\geq 0}$ in $\mathbb{C}\{\{x-a\}\}$

**Proposition 2.2.** (a) Assume that $f \in \mathbb{C}\{\{x\}\}_\mathbb{P}$ has radius of convergence $r > 0$ and that $a \in \mathbb{C}$ with $|a| < r$.

Then for all $T \in \mathbb{P}$ the sequence $(\langle g_n, (x-a)^T \rangle)_{n \geq 0}$ is convergent in $\mathbb{C}$. If its limit is denoted by $\gamma_T(a)$, then
$$f := \sum_{T \in \mathbb{P}} \gamma_T(a) \cdot (x-a)^T \in \mathbb{C}\{\{x-a\}\}_\mathbb{P}$$
It is called the expansion of $f$ around the place $a$.

Moreover $\gamma_T(a) = \langle g, (x-a)^T \rangle =$
$$= \sum_{U \in \mathbb{P}} \langle f, x^U \rangle \, (U/T) \cdot a^{(\deg(U) - \deg(T))}$$
and the expansion of $f$ around $a$ has radius of convergence $\geq r - |a|$.

Special cases.
If $T = \mathbf{1} =$ empty tree,
then $(U/T) = 1$ for all $n \in \mathbb{P}$ and $\gamma_\mathbf{1}(a) = \sum_{U \in \mathbb{P}} \gamma_U \cdot a^{(\deg(U))} = f(a)$

If $T = x$, then $(U/x) = \deg(U)$, if $n \neq \mathbf{1}$ and
$$\gamma_x(a) = \sum_{\deg(U) \geq 1} \gamma_U \cdot \deg(U) \cdot a^{(\deg(U) - 1)}$$

# 3 The notion of planar analytic functions.

Let $D$ be a region in the complex plane $\mathbb{C}$ and $F = (f_a)_{a \in D}$ a collection of planar power series $f_a$ in the variable $x - a$ for all $a \in G$.
Definition: F is called planar analytic function on D if
(i) $f_a$ is a convergent planar power series which means that the radius of convergence of $f_a$ is $> 0$
(ii) if $a, b \in G$ and if the radius of convergence rad $(f_a) > |b-a|$, then the expansion of $f_a$ around $b$ is equal to $f_b$.

Also $f_a$ is called the germ of F at a. There is a sheaf $\vartheta_{\mathbb{P},D}$ on D for which the stalk in $a \in D$ consists of the convergent planar power series in the variable $(x-a)$.

The sections in $\Gamma(D, \vartheta_{\mathbb{P},D})$ are the planar analytic functions on D.

**Proposition (3.1).** Let $f \in \mathbb{C}\{\{x\}\}_\mathbb{P}$, $r = \text{rad}(f) > 0$ and $G = \{a \in \mathbb{C} : |a| < r\}$.
Let $f_a$ be the expansion of $f$ around $a$ and $F = (f_a)\, a \in G$.
Then F ist a planar analytic function on G.

To prove this statement we have to show the following: Assume that $a, b \in \mathbb{C}$ with $|a| < R$, $|b| < R$ where R is the radius of convergence of a planar series $f$ in x. Assume more over that $|b-a|$ is smaller than the radius of convergence of the expansion $f_a$ of $f$ around a. Then the expansion $f_b$ of $f$ around $b$ coincides with the eypansion of $f_a$ around b. It means that the expansions do not depend on the "paths" within the circle of convergence.

Thus we have to show that

$$\gamma_T(b) := \sum_{U \in \mathbb{P}} \gamma_U(U/T)\, b^{(\deg(U) - \deg(T))}$$

is equal to

$$\delta_T(b) := \sum_{S \in \mathbb{P}} \gamma_U(a)\, (S/T)\, (b-a)^{(\deg(S) - \deg(T))}$$

Substituting $\gamma_T(a)$ into this expression gives

$$\delta_T(b) = \sum_{S \in \mathbb{P}} \sum_{S \in \mathbb{P}} \gamma_U(U/T)\, (S/U)\, a^{\deg(U) - \deg(T)}\, (b-a)^{(\deg(S) - \deg(T))}$$

The crucial observation in proving $\gamma_T(b) = \delta_T(b)$ is the following relation:

Let $S, T \in \mathbb{P}$, $\deg(S) = n$, $\deg(T) = m$, and $a, b \in \mathbb{C}$. Then

$$\sum_{U \in \mathbb{P}} (S/U)(U/T)\, a^{n - \deg(U)} \cdot (b-a)^{\deg(U) - m} = (S/T)\, b^{n-m}$$

The left hand side is equal to

$$\sum_{k=0}^{n-m} \left( \sum_{\deg(U) = m+k} (S/U)(U/T)\, a^{n-m-k}\, (b-a)^k \right)$$

and

$$\sum_{\deg(U) = m+k} (S/U)(U/T) = (S/T)\, (n-m/k)$$

The left hand side of this expression
is $\# M$ with $M := \{(J, I) : J, I \leq L(S), \# I = m, \# J = m + k,$ the contraction $S|I$ is equal to $T\}$ which is equal to $(S/T)\,(n-m/k)$ because for fixed $I$ with $S|I = T$ there are $(n-m/k)$ subsets $J$ with $(J, I) \in M$.

It follows that the left hand side is equal $M$

$$\sum_{k=0}^{n-m} (S/T)\,(n-m/k)\, a^{n-m-k}\, (b-a)^k$$

witch is equal to $(S/T)\, b^{n-m.}$

Let $F = (f_a)_{a \in D}$ be a planar analytic function on $D$.
For $T \in \mathbb{P}$ the function $h_T : D \to \mathbb{C}$ given by $h_T(a) := \langle f_a, (x-a)^T \rangle$ is analytic, it is called the $T$ – th coefficient $\langle F, T \rangle$ of $F$.

## 4 Exponential functions

Let $\exp_k(x) \in \mathbb{C}\{\{x\}\}_\mathbb{P}$ be the k –ary exponential series where $k \in \mathbb{N}_k \geq 2$, see [G1]. Then from the recursion formula for the coefficients $A_T(k)$ of $\exp_k(x)$, see [G1], p. 352, it follows that all coefficients $A_T(k)$ of $\exp_k(x)$ are real and positive and

$$\sum_{T \in \mathbb{P}(n)} A_{kT} = 1/n$$

where $\mathbb{P}(n)$ is the set of trees in $\mathbb{P}$ of degree n.
It follows that the radius of convergence of $\exp_k(x)$ is $\infty$.

**Remak 4.1.** It is not obvious that the radius of convergence
of $\log_k(1+x)$ is equal to 1.

Let $\lambda \in \mathbb{C}$. Then $f = \exp_k(\lambda+x)$ is a planar power series in $\mathbb{C}\{\{x+\lambda\}\}_\mathbb{P}$ which has infinite radius of convergence. Thus f has expansion around 0,
and $\quad f(\lambda+x) = \sum \alpha_{k,T}(\lambda) \, x^T$

with $\quad \alpha_{k,T}(\lambda) = \sum_{U \in \mathbb{P}} (U/T) \, \alpha_{k,n} \cdot \lambda^{(\deg(U) - \deg(T))}$

**Proposition 4.2.** $\exp_k(\lambda+x) = e^\lambda \cdot \exp_k(x)$

Proof: Proceed as in the proof of Proposition 5.1 in [G5]

Work in $A = \mathbb{C}[\lambda] \widehat{\otimes} \, \mathbb{C}\{\{x\}\}_\mathbb{P}$

Let $A_n := \mathbb{C}$ - vector space generated by $\{\lambda^i \cdot x^T : i + \deg(T) = n\}$
Each $f \in A$ has a unique expansion

$$f = \sum_{n=0}^{\infty} f_n,$$

with $f_n \in A_n$.

There is a unique $f \in A$ such that $f = 1 + x + \sum_{n=2}^{\infty} f_n \quad f_n \in A$

$$f^k(\lambda, x) = f(k\cdot\lambda, k\cdot x)$$

One can show that both $e^\lambda \cdot \exp_k(x)$ and $\exp_k(\lambda+x)$ do satisfy both conditions.
This shows that they are equal.

**Corollary 4.3.** Let $n = \deg(T)$. Then

$$(\alpha_T/m!) = \sum_{U \in \mathbb{P}_{n+m}} (U/T) \, \alpha_U$$

**Example 4.4.** $T = x \cdot x^2$, $m = 1$:

$\mathbb{P}_4 = \{x \cdot (x \cdot x^2), x(x^2 x), x^2 \cdot x^2, (x \cdot x^2) x, (x^2 x) \cdot x\} = \{U_1,..., U_5\}$
and $(U_1/T) = 4, (U_2/T) = 3, (U_3/T) = 2, (U_4/T) = 1, (U_5/T) = 0$

Thus $\alpha_T / 1! = 4 \alpha_{U_1} + 3 \alpha_{U_2} + 2 \alpha_{U_3} + \alpha_{U_4}$

$\alpha_{ni} = 1/4! \cdot 1/7$ if $i \neq 3$ and $\alpha_{n3} = \frac{1}{4}! \, 1/7 \cdot 3$

Thus $\alpha_T = 4 \alpha_{U_1} + 3 \alpha_{U_2} + 2_{U_3} + \alpha_{U_4}$

$$= \frac{1}{4! \cdot 7} (4 + 3 + 2 + 2 \cdot 3 + 1) = \frac{14}{4! \cdot 7} = \frac{2}{4!} = \frac{1}{3!} \cdot \frac{1}{2}$$

# 5 Planar zeta function and Gamma function

Fix $k \in \mathbb{N}_{\geq 2}$ and let $\exp(s) = \exp_k(s)$ be the planar exponential series of arity $k$ in the planar variable $s$.

For an integer $n \in \mathbb{N}_{\geq 1}$ let $n^s := \exp(s \cdot \log(n))$.

If $\exp(s) = \sum_{T \in \mathbb{P}} \alpha_T \, s^T$, one gets

$$n^s = \sum_{T \in \mathbb{P}} \alpha_T \cdot (\log(n))^{\deg(T)} \cdot s^T$$

Let $r \in \mathbb{C}$. Then $n^{s+r} = n^s \cdot n^r$

**Proposition 5.1.** There is a planar analytic function $\zeta_\mathbb{P}$ on $D = \mathbb{C} - \{\lambda \in \mathbb{R} : \lambda \leq 1\}$ such that for $r \in D$, $\mathrm{Re}(r) > 1$, one has the expansion $\zeta_\mathbb{P} = \sum_{n=1}^{\infty} n^{-r} (n^{-(s-r)})$ around $r$.

More precisely

$$\zeta_\mathbb{P} = \sum_{n=1}^{\infty} \sum_{T \in \mathbb{P}} 1/n^r \, \alpha_T \, (-\log(n))^{\deg(T)} \, (s-r)^T$$

$$\zeta_\mathbb{P} = \sum_{T \in \mathbb{P}} \beta_{T,r} \cdot (s-r)^T$$

with $\beta_{T,r} = \left( \sum_{n=1}^{\infty} 1/n^r \, (-\log(n))^{\deg(T)} \right) \cdot \alpha_T$

The radius of convergence of the expansion of $\zeta_\mathbb{P}$ around $r$ is $(r-1)$.

Open question: what is the relation between $\zeta_\mathbb{P}(s)$ and $\zeta_\mathbb{P}(1-s)$ ?

The definition of the Gamma function $\Gamma(s)$ of Euler by the integral

$$\int_0^\infty e^t \cdot t^{s-1} \, dt$$

has an immediate extension as a planar analytic function in the

domain $\{a \in \mathbb{C} \cdot \mathrm{Re}(a) > 0\}$. This can be done as follows:

Let $r \in \mathbb{C}$, $\mathrm{Re}(r) > 0$. We work in the algebras $\mathbb{C}[|t|] \otimes \mathbb{C}\{\{s-r\}\}$.

Then $t^{s-1} = t^{(s-r)} \cdot t^{(r-1)}$ where $t^{(s-r)}$ is the planar function $\exp((s-r)\log t)$ and where $\exp$ is the $k$-ary planar exponential series for some $k \in \mathbb{N}, k \geq 2$.

Thus

$$t^{s-1} = \sum_{T \in \mathbb{P}} t^{(r-1)} \, \alpha_T \cdot (\log t)^{\deg(T)} \cdot (s-r)^T$$

Now for fixed $T \in \mathbb{P}$ the integral.

$$\beta_T(r) := \int_0^\infty e^{-t} t^{(r-1)} \cdot \alpha_T (\log t)^{\deg(T)} dt$$

is well defined and $\Gamma_\mathbb{P}(s)$ is defined to have the expansion around $r$ given by the planar series

$$\sum_{T \in \mathbb{P}} \beta_T(r) (s-r)^T$$

The open problem is again what kind of functional equation holds for $\Gamma_\mathbb{P}$.

Certainly $\Gamma_\mathbb{P}(s+1)$ is different from $s \cdot \Gamma(s)$.

# 6  The multiplicative inverse $F^{-1}$

Let $F = (f_a)_{a \in D}$ be a planar analytic function on D and assume that $f_a(a) \neq 0$ for all $a \in D$. Then the left inverse $f_a^{-1}$ of $f_a$ exists in $\mathbb{C}\{\{x-a\}\}_\mathbb{P}$ which satisfies

$$f_a^{-1} \cdot f_a = \underline{1}$$

Let $F^{-1} = (f_a^{-1})_{a \in D}$

**Proposition 6.1.**  $F^{-1}$ is a planar analytic function on D.

The same statement holdes for $^{-1}F = (^{-1}f_a)_{a \in D}$ where $^{-1}f_a$ is the right inverse of $f_a$.

# The planar function $(x^{-1})$

Let $D = \mathbb{C} - \{0\}$, $a \in D$ and $f_a = ((x-a) + a)^{-1}$ be the left inverse of $(x-a) + a$

in $\mathbb{C}\{\{x-a\}\}_\mathbb{P}$ which is $a^{-1}(1 + x-a/a)^{-1} = \sum_{n=0}^\infty 1/a^n \, (x-a)^{n)}$

where $(x-a)^{(n+1))} := (x-a)^{n)} \cdot (x-a) = (x-a)^{C_{n+1}}$

where $C_{n+1}$ is the right comb of degree $n + 1$.

Then $(f_a) \, a + D$ is a planar analytic function also denoted by $(x^{-1})$

# 8 The planar square root function ($_\mathbb{P}\sqrt{x}$)

Let $D = \mathbb{C} - \{\lambda \in \mathbb{R} : \lambda \leq 0\}$ and $\sqrt{\phantom{z}} : D \to \mathbb{C}$ the analytic function with $\sqrt{1} = 1$ and $(\sqrt{z})^2 = z$ for all $z \in D$.

For all $a \in D$ there is a unique planar series $f_a \in \mathbb{C}\{\{x-a\}\}_\mathbb{P}$ such that the constant term of $f_a$ is $\sqrt{a}$ and $f_a^2 = a + (z-a)$. One can show that the radius of convergence of $f_a$ is equal to $|a|$. The planar analytic function $_\mathbb{P}\sqrt{x} := (f_a)_{a \in D}$ is called the planar square root function.

Let $\quad g = \sum_{\substack{T \in \mathbb{P}' \\ \deg(T) \geq 1}} x^T$

where $\mathbb{P}'$ ist he set of trees T in $\mathbb{P}$ for each vertex has arity less or equal to 2.

It is easy to show that $g^2 = g - x$ and thus $(g - \tfrac{1}{2})^2 = 1/4 - x$

It follows that $(1-2g)^2 = 1 - 4x$

Let $h = g(-1/4\, x)$. Then $(1-2h)^2 = 1 + x$ and $h = \sum_{\substack{T \in \mathbb{P} \\ \deg(T) \geq 1}} (-1/4)^{\deg T} x^T$

The radius of convergence of g is equal to 1/4 because all coefficients of g are positive and the classical series $[[g]] = \tfrac{1}{2} - \tfrac{1}{2}\sqrt{1-4x}$ and $\sqrt{1-4x}$

$$= \sum_{n=0}^{\infty} (\tfrac{1}{2}/n)(-4x)^n$$

where $(\tfrac{1}{2}/n)$ is the binominal coefficient of ½ over n, which has ¼ as radius of convergence.